\documentclass[a4paper, 12pt]{article}
\usepackage{amssymb}
\usepackage{amstext}
\usepackage{amsmath}
\usepackage{amscd}
\usepackage{latexsym}
\usepackage{theorem}
\usepackage{color}

\theoremheaderfont{\scshape}
\theorembodyfont{\itshape}
\newtheorem{thm}{\bfseries Theorem}[section]
\newtheorem{prop}[thm]{\bfseries Proposition}
\newtheorem{lemma}[thm]{\bfseries Lemma}
\newtheorem{cor}[thm]{\bfseries Corollary}
\theorembodyfont{\rmfamily}
\newtheorem{defn}[thm]{\bfseries Definition}
\newtheorem{ex}[thm]{\bfseries Example}
\newtheorem{rem}[thm]{\bfseries Remark}

\newtheorem{pf}{Proof.}

\numberwithin{equation}{section}

\def\N{\mathbb{N}}

\def\mod{\mathrm{mod}}
\def\CM{\mathrm{CM}}
\def\sCM{{\underline{\mathrm{CM}}}}
\def\Hom{\mathrm{Hom}}
\def\sHom{{\underline{\mathrm{Hom}}}}

\def\Ext{\mathrm{Ext}}

\def\End{\mathrm{End}}
\def\sEnd{{\underline{\mathrm{End}}}}

\def\coker{\mathrm{Coker}}
\def\ker{\mathrm{Ker}}
\def\image{\mathrm{Im}}

\def\rank{\mathrm{rank}}
\def\m{\mathfrak m}

\def\p{\mathfrak p}

\def\P{\mathfrak P}

\def\length{\mathrm{length}}

\def\Spec{\mathrm{Spec}}

\def\F{{\mathcal F}}
\def\qed{$\Box$}

\def\phipsi{\fontsize{6pt}{6pt} {\binom{\phi}{\psi}}}
\def\01{\fontsize{6pt}{6pt} {\binom{0}{1}}}

\def\u#1{\underline{#1}}

\begin{document}


\title{\bf Stable degenerations \\ 
of Cohen-Macaulay modules} 

\author{
{Yuji Yoshino}\\
{\small Department of Mathematics}\vspace{-1mm}\\ 
{\small Faculty of Science}\vspace{-1mm}\\   
{\small Okayama University}\vspace{-1mm}\\ 
{\small Okayama 700-8530, Japan}\vspace{-1mm}\\ 
{\small e-mail: {\tt yoshino@math.okayama-u.ac.jp}}\vspace{1mm}\\
}
\date{\empty}
\maketitle

\medskip

\begin{abstract} 
As a stable analogue of degenerations, we introduce the notion of stable degenerations for Cohen-Macaulay modules over a Gorenstein local algebra. 
We shall give several necessary and/or sufficient conditions for the stable degeneration. 
These conditions will be helpful to see when a Cohen-Macaulay module degenerates to another. 
\end{abstract}

\thispagestyle{empty}

\section{Introduction}\label{Intro}

Let  $k$  be an algebraically closed field and let  $R$  be a finite-dimensional associative $k$-algebra with basis $e_1=1, e_2, \ldots , e_d$. 
Then the structure constants $c_{ijk}$  are defined by 
$e_i \cdot e_j = {\sum}_k c_{ijk} e_k$. 
The module variety $\mathrm{Mod} _R^n$  of $n$-dimensional left $R$-modules consists of the $d$-tuple  $x = (x_1, x_2 ,\ldots, x_d)$ of $n \times n$-matrices with entries in  $k$  such that  $x_1$  is the identity and 
$x_ix_j = \sum _k c_{ijk} x_k$  holds for all $i,j$.  
It is easy to see that  $\mathrm{Mod}_R^n$  is an affine variety and the general linear group $\mathrm{GL}_n(k)$  acts on  $\mathrm{Mod} _R^n$  by simultaneous conjugation. 
In such a situation each $\mathrm{GL}_n(k)$-orbit corresponds to an isomorphism class of $n$-dimensional modules. 
We denote by  $\mathcal{O}(M)$  the $\mathrm{GL}_n(k)$-orbit of the point in $\mathrm{Mod}_R^n$  corresponding to an $R$-module $M$. 
Then we say that  $M$  degenerates to $N$  if  $\mathcal{O}(N)$  is contained in the Zariski closure of  $\mathcal{O}(M)$. 
This is the definition of degeneration of modules over a finite dimensional algebra.
See \cite{B}, \cite{R} and \cite{Z}. 

However, if we want to consider the degeneration of modules over a noetherian algebra $R$ which is not necessarily  finite-dimensional over $k$, then this definition is not applicable, since each module is not necessarily a finite-dimensional $k$-vector space any more even if it is assumed to be finitely generated over $R$,  and hence there is no good way to define appropriate module varieties. 
By this reason we have proposed a scheme-theoretical definition of degeneration in the previous paper \cite{Y3}. 
Our definition is as follows: 

Let  $R$  be an associative $k$-algebra where $k$ is any field. 
We set  $V = k[t]_{(t)}$ and  $K=k(t)$  with  $t$  being a variable over  $R$ hence over  $k$. 
We denote by  $\mod (R)$  the category of all finitely generated left $R$-modules and $R$-homomorphisms. 
Then we have the natural functors 
$$
\begin{CD}
\mod (R)  @<{r}<< \mod (R \otimes _k V) @>{\ell}>> \mod (R \otimes _k K), 
\end{CD}
$$ 
where  $r = - \otimes _V V/tV$  and  $\ell = - \otimes _V K$. 
Then for modules  $M, N \in \mod (R)$, we say that  $M$ degenerates to $N$ if there is a module  $Q \in \mod (R \otimes _k V)$  that is $V$-flat such that  $\ell (Q) \cong M \otimes _k K$  and  $r(Q) \cong N$. 
The module  $Q$, regarded as a bimodule  ${}_R Q_V$, is a flat family of $R$-modules with parameter in  $V$. 
At the closed point in the parameter space $\Spec V$, the fiber of $Q$ is $N$, which is a meaning of the isomorphism  $r(Q) \cong N$. 
On the other hand, the isomorphism $\ell (Q) \cong M \otimes _k K$  means that   the generic fiber of $Q$  is essentially given by $M$. 
In \cite{Y3} we have shown that the latter definition of degeneration agrees with the former if  $R$  is a finite-dimensional $k$-algebra.

In the present paper we want to consider the stable analogue of degenerations for Cohen-Macaulay modules over a Gorenstein local rings. 

Let  $R$  be a Gorenstein commutative  local $k$-algebra where $k$ is any field, and set  $V = k[t]_{(t)}$ and  $K=k(t)$  as above. 
We denote by  $\CM (R)$  the category of all (maximal) Cohen-Macaulay  $R$-modules and all $R$-homomorphisms, and denote by $\sCM(R)$  the stable category of  $\CM(R)$. 
Note that  $\sCM (R)$  has a structure of a triangulated category. 
Then, similarly to the above, there are triangle functors 
$$
\begin{CD}
\sCM(R)  @<{\mathcal{R}}<< \sCM (R \otimes _k V) @>{\mathcal{L}}>> \sCM (R \otimes _k K), 
\end{CD}
$$ 
where  $\mathcal{R}$  and  $\mathcal{L}$  are induced respectively by $- \otimes _V V/tV$  and  $ - \otimes _V K$. 
Then we define that  $\u{M} \in \sCM (R)$ stably degenerates to $\u{N} \in \sCM (R)$ if there is a  Cohen-Macaulay module  $\u{Q} \in \sCM (R \otimes _k V)$  such that  $\mathcal{L} (\u{Q}) \cong \u{M \otimes _k K}$  and  $\mathcal{R}(\u{Q}) \cong \u{N}$. 
See Definition \ref{def stable deg} for more detail. 

The aim of this paper is to clarify the meaning of this definition by showing several implications. 
To explain this, let  $(R, \m, k)$  be a Gorenstein complete local $k$-algebra and assume  for simplicity that  $k$  is an infinite field. 
For Cohen-Macaulay $R$-modules $M$ and  $N$ we consider the following four conditions: 
\begin{itemize}
\item[$(1)$]
$R ^m \oplus M$  degenerates to  $R^n \oplus N$  for some  $m, n \in {\Bbb N}$. 
\item[$(2)$]
There is a triangle  
$\u{Z} \overset{\binom{\u{\phi}}{\u{\psi}}}\to \u{M} \oplus \u{Z} \to \u{N} \to \u{Z}[1]$  in  $\u{\rm CM}(R)$,  where $\u{\psi}$  is  a nilpotent element of  $\u{\rm End} _R (Z)$. 

\item[$(3)$]
$\u{M}$ stably degenerates to  $\u{N}$.

\item[$(4)$]
There exists an  $X \in \CM (R)$  such that  $M \oplus R^m \oplus X$  degenerates to $N \oplus R^n \oplus X$ for some  $m, n \in \N$.

\end{itemize}

In this paper, we shall prove the following implications and equivalences of these conditions:

\begin{itemize}
\item
In general, $(1) \Rightarrow (2) \Rightarrow (3) \Rightarrow (4)$. 
(Theorems \ref{implication} and \ref{oplusX}).  

\item
If  $R$  is an artinian ring, then $(1)\Leftrightarrow(2)\Leftrightarrow(3)$. 
(Theorem \ref{implication}).  

\item
If  $R$  is an isolated singularity, then  $(2) \Leftrightarrow (3)$. 
(Theorem \ref{isolated singularity}). 
\end{itemize}
However there is a counter example to the implication $(2) \Rightarrow (1)$  in the case when $R$  is an isolated singularity of Krull dimension one (Example \ref{ex 2->1}).  
We also give an example of a Gorenstein ring of Krull dimension zero for which the implication $(4) \Rightarrow (3)$ fails (Example \ref{Riedtmann}). 
Furthermore, as a consequence of the implication $(3) \Rightarrow (4)$, we can show that the stable degeneration gives rise to a well-defined  partial order on the set of isomorphism classes of  Cohen-Macaulay modules (Theorem \ref{st deg order}).

After giving some preliminary consideration for degenerations in Section 2, we make several remarks on the ring $R \otimes _k K$ in Section 3. 
The remaining part of the paper is devoted to giving the proofs of the results mentioned above and to constructing the examples.

\section{Review of degenerations}\label{review}

Let us recall the precise definition of degeneration of finitely generated modules over a noetherian algebra, which is given in our previous paper \cite[Definition 2.1]{Y3}.

\begin{defn}\label{defDVR} 
Let  $R$  be a noetherian algebra over a field  $k$, and let  $M$ and $N$  be finitely generated $R$-modules. 
We say that $M$ degenerates to  $N$, or  $N$  is a degeneration of  $M$, if 
there is a discrete valuation ring  $(V, tV, k)$  that is a $k$-algebra (where $t$  is a prime element) and a finitely generated  left  $R\otimes _k V$-module  $Q$  which satisfies the following conditions: 
\begin{itemize}
\item[(1)]
$Q$  is flat as a $V$-module.
\item[(2)]
$Q/tQ \cong N$  as a left $R$-module.
\item[(3)]
$Q[\frac{1}{t}] \cong M \otimes _k V[\frac{1}{t}]$  as a left $R\otimes _k V[\frac{1}{t}]$-module.
\end{itemize}
\end{defn}

In the previous paper(\cite[Theorem 2.2]{Y3}) we have proved the following theorem.

\begin{thm}\label{deg}
The following conditions are equivalent for finitely generated left $R$-modules $M$ and  $N$.
\begin{itemize}
\item[$(1)$]
$M$ degenerates to $N$.
\item[$(2)$]
There is a short exact sequence of finitely generated left $R$-modules
$$
0 \to Z \overset{\phipsi}\longrightarrow M \oplus Z \to N \to 0,
$$
such that the endomorphism  $\psi$  of  $Z$ is nilpotent, i.e. $\psi ^n = 0$  for  $n \gg 1$.
\end{itemize}
\end{thm}

By virtue of this theorem together with a theorem of Zwara \cite[Theorem 1]{Z}, we see that if  $R$  is a finite-dimensional algebra over  $k$, then our definition of degeneration agrees with the original one using module varieties which is mentioned in the beginning of Introduction. 
We also remark from this theorem that we can always take  $k[t]_{(t)}$  as $V$ in Definition \ref{defDVR}. 
(See \cite[Corollary 2.4]{Y3}.)

In the rest of the paper we mainly treat the case when $R$  is a commutative ring.

\begin{rem}\label{if degerates}
Let  $R$  be a commutative noetherian algebra over $k$, and suppose that a finitely generated $R$-module  $M$ degenerates to a finitely generated $R$-module  $N$. 
Then the following hold. 

\begin{itemize}
\item[(1)]
The modules $M$ and $N$ give the same class in the Grothendieck group, i.e. 
$[M] = [N]$  as elements of  $K_0 (\mod (R))$. 
This is actually a direct consequence of Theorem \ref{deg}. 
In particular, $\rank \  M = \rank \ N$ if the ranks are defined for $R$-modules. 
Furthermore, if  $(R, \m)$  is a local ring, then $e(I, M)= e(I, N)$  for any $\m$-primary ideal $I$, where  $e(I, M)$  denotes the multiplicity of $M$ along $I$. 

\item[(2)] 
If  $L$  is an $R$-module of finite length, then we have the following inequalities of lengths for any integer  $i$:  
$$
\begin{cases}
&\length _R(\Ext _R ^i(L, M) ) \leqq \length _R(\Ext _R ^i (L, N)), \\ 
&\length _R(\Ext _R ^i(M, L) ) \leqq \length _R(\Ext _R ^i (N, L)). 
\end{cases}
$$
See \cite[Lemma 4.5]{Y2}.
In particular, when $R$  is a local ring, 
then  
$$
\nu (M) \leqq \nu (N), \quad  \beta _i (M) \leqq \beta _i(N) \ \ \text{and}\ \  \mu ^i (M) \leqq \mu ^i (N) \ \ (i\geqq 0), 
$$
where  $\nu$, $\beta _i$ and $\mu ^i$ denote the minimal number of generators, the $i$th Betti number and  the $i$th Bass number respectively. 
\end{itemize}
\end{rem}


Adding to these remarks we can prove the following result concerning Fitting ideals. 
Before stating the theorem, we recall the definition of the Fitting ideal of a finitely presented module. 
Suppose that a module  $M$  over a commutative ring $A$  is given by a finitely free presentation 
$$
\begin{CD}
A^{m} @>{C}>> A^{n} @>>> M @>>> 0,  
\end{CD}
$$
where  $C$  is an $n \times m$-matrix with entries in  $A$. 
Then recall that the $i$th Fitting ideal  $\F_i^A (M)$  of $M$  is defined to be the ideal $I_{n-i}(C)$  of  $A$  generated by all the $(n-i)$-minors of the matrix $C$.
(We use the convention that  $I_r(C) =A$  for  $r \leqq 0$ and $I_r(C) = 0$ for  $r> \min \{m, n\}$.) 
It is known that $\F _i^A(M)$  depends only on $M$ and  $i$, and independent of the choice of free presentation. 
The following lemma will be used to prove the theorem. 

\begin{lemma}\label{Fitting}
Let  $f: A \to B$ be a ring homomorphism and let  $M$  be an $A$-module which possesses a finitely free presentation. 
Then $\F _i ^B (M \otimes _A B) = f(\F_i ^A (M))B$  for all $i \geqq 0$. 
\end{lemma}

\begin{pf}
If $M$ has a presentation 
$$
\begin{CD}
A^{m} @>{(c_{ij})}>> A^{n} @>>> M @>>> 0,  
\end{CD}
$$
then $M \otimes _A B$  has a presentation 
$$
\begin{CD}
B^{m} @>{(f(c_{ij}))}>> B^{n} @>>> M \otimes _A B @>>> 0.  
\end{CD}
$$
Thus $\F _i ^B (M \otimes _A B) = I_{n-i}(f(c_{ij}))= f(I_{n-i}(c_{ij}))B = f(\F_i ^A (M))B$. 
\qed
\end{pf}

\begin{thm}\label{decreaseFitting}
Let  $R$  be a noetherian commutative algebra over $k$, and $M$ and $N$ finitely generated $R$-modules. 
Suppose  $M$  degenerates to $N$. 
Then we have  $\F_i^R (M)  \supseteqq \F_i^R (N)$  for all  $i \geqq 0$. 
\end{thm}

\begin{pf}
By the assumption there is a finitely generated $R \otimes _kV$-module $Q$  such that $Q_t \cong  M \otimes _kK$ and  $Q/tQ \cong N$, where  $V = k[t]_{(t)}$ and  $K = k(t)$. 
Note that  $R \otimes _k V \cong S^{-1}R[t]$  where  $S = k[t]\backslash (t)$. 
Since  $Q$  is finitely generated, we can find a finitely generated  $R[t]$-module  $Q'$  such that  $Q' \otimes _{R[t]} (R \otimes _k V) \cong Q$. 
For a fixed integer  $i$ we now consider the Fitting ideal  $J :=\F_i^{R[t]} (Q') \subseteqq R[t]$. 
Apply Lemma \ref{Fitting} to the ring homomorphism  $R[t] \to R =R[t]/tR[t]$, and noting that  $Q' \otimes _{R[t]} R \cong N$,  we have  
\begin{equation}\label{J}
\F_i^R (N) = J + tR[t] /tR[t] \hphantom{CCCC} 
\end{equation}
as an ideal of  $R = R[t] /tR[t]$. 
On the other hand, applying Lemma \ref{Fitting} to  $R[t] \to R\otimes _k K = T^{-1}R[t]$ where  $T = k[t] \backslash \{ 0\}$, we have $\F_i ^R (M) T^{-1}R[t]   =  JT^{-1}R[t]$. 
Therefore there is an element  $f(t) \in T$  such that $f(t)J \subseteq \F_i^R(M)R[t]$. 
  
Now to prove the inclusion  $\F_i^R (N) \subseteqq \F _i^R (M)$, take an arbitrary element  $a \in \F_i^R (N)$. 
It follows from (\ref{J}) that there is a polynomial of the form 
$a + b_1 t + b_2 t^2 + \cdots + b_r t^r\ (b_i \in R)$ that belongs to $J$. 
Then, 
we have  $f(t)(a + b_1 t + b_2 t^2 + \cdots + b_r t^r) \in \F_i ^R (M) R[t]$. 
Since  $f(t)$  is a non-zero polynomial whose coefficients are all in $k$,    
looking at the coefficient of the non-zero term of the least degree in the polynomial $f(t)(a + b_1 t + \cdots + b_r t^r)$,  we have that  $a \in \F_i^R (M)$. 
\qed
\end{pf}

\section{Remarks for the rings $R\otimes _k K$.}

In this paper we are interested in the stable analogue of degenerations of  Cohen-Macaulay modules over a commutative Gorenstein local ring. 
For this purpose, all rings considered in the rest of the paper are commutative. Furthermore  $(R, \m, k)$  always denotes a Gorenstein local ring which is a $k$-algebra, and  $V = k[t]_{(t)}$  and  $K = k(t)$ where  $t$ is a variable.

We note that $R \otimes _k V$  and  $R \otimes _k K$  are not necessarily local rings but they are Gorenstein rings.

\begin{rem}\label{remark on RotimesK}
Let  $S = \{ f(t) \in k[t] \ \vert \ f(0) \not= 0\}=k[t]\backslash (t)$ and 
$T = k[t] \backslash \{0\}$, which are multiplicatively closed subsets of  $k[t]$, hence of  $R[t]$. 
Then, by definition, we have 
$$
R\otimes _k V = S^{-1} R[t], \quad R\otimes _k K = T^{-1} R[t]. 
$$
In particular, both rings are noetherian. 
Furthermore, there are natural mappings $R \to R\otimes _k V$  and $R \to R\otimes _k K$, which are flat ring homomorphisms. 
For any $\p \in \Spec (R)$, the fibers of these homomorphisms are
$$
 \kappa (\p) \otimes _R (R\otimes _k V) = S^{-1} \kappa (\p)[t], \quad 
 \kappa (\p) \otimes _R (R\otimes _k K) = T^{-1} \kappa (\p)[t], 
$$ 
which are regular rings. 
Hence we see that  $R\otimes _k V$ and $R\otimes _k K$  are Gorenstein as well as  $R$. 
At the same time,  we have the equality of Krull dimension, 
$$
\dim R\otimes _k V = \dim R +1, \quad \dim R\otimes _k K = \dim R.
$$

If $\dim R =0$ (i.e. $R$  is artinian), then the rings  $R \otimes _k V$  and  $R \otimes _k K$  are local. 
In fact, the ideal  $\m (R \otimes _kV)$  of  $R \otimes _k V$ is nilpotent, and $(R\otimes _k V)/\m(R\otimes _k V) \cong V$, hence $(\m, t)(R\otimes _kV)$ is a unique maximal ideal of  $R \otimes _k V$. 
By the same reason,  $\m(R \otimes _k K)$ is a unique maximal ideal of  $R \otimes _k K$. 

However we should note that  $R \otimes _k V$  and  $R \otimes _k K$ will never be local rings if $\dim R > 0$.
Actually, if there is a prime ideal  $\p$  with $\dim R/\p =1$, then taking an  $x \in R$ so that  $x \not\in \p$, we have a maximal ideal  $(\p, xt -1) R \otimes _kV$ (resp.  $(\p, xt -1) R \otimes _kK$), which is distinct from the maximal ideal  $(\m, t)R \otimes _k V$ (resp. $\m (R \otimes _kK)$).   
\end{rem}

Since  $R \otimes _k K$  is non-local, there may be a lot of projective modules which are not free. 
The following example gives such one of them.

\begin{ex}\label{nonfree-projective} 
Let  $R = k[[ x, y]]/(x^3-y^2)$. 
It is known that the maximal ideal $\m = (x, y)$ is a unique non-free indecomposable Cohen-Macaulay module over  $R$. 
See \cite[Proposition 5.11]{Y1}.
In fact it is given by a matrix factorization of the polynomial  $x^3-y^2$: 
$$
(\varphi, \psi) = 
\left(
\begin{pmatrix}
y   & x \\
x^2 & y \\
\end{pmatrix}, \quad 
\begin{pmatrix}
y   & -x \\
-x^2 & y \\
\end{pmatrix}
\right). 
$$
Therefore there is an exact sequence  
$$
\begin{CD}
R^2 @>{\varphi}>> R^2 @>>> \m @>>> 0. \\
\end{CD}
$$
Now we deform these matrices and consider the pair of matrices over  $R \otimes _k K$: 
$$
(\Phi, \Psi) = 
\left(
\begin{pmatrix}
y -xt   & x - t^2\\
x^2 & y + xt \\
\end{pmatrix}, \quad 
\begin{pmatrix}
y +xt  & -x + t^2\\
-x^2 & y -xt \\
\end{pmatrix}
\right). 
$$
Define the $R \otimes _k K$-module $P$ by the following exact sequence:
$$
\begin{CD}
(R \otimes _k K)^2 @>{\Phi}>> (R \otimes _k K)^2 @>>> P @>>> 0. \\
\end{CD}
$$

\vspace{6pt}
\noindent
We claim that 
\underline{$P$ is a projective module of rank one  over  $R \otimes _k K$ but non-free}.
(Hence the Picard group of  $R \otimes _k K$  is non-trivial.)
\vspace{6pt}

To prove that  $P$ is projective of rank one. let  $\P$ be a prime ideal of  $R \otimes _k K$. 
If  $x \not\in \P$, then  $x^2$  is a unit in $(R \otimes _k K)_{\P}$, therefore the matrix  $\Phi _{\P}$  over  $(R \otimes _k K)_{\P}$  is equivalent to  
$\begin{pmatrix} 0& 0 \\ 1 & 0 \end{pmatrix}$ under elementary transformations, hence  $P_{\P} = \coker (\Phi _{\P}) \cong \coker (\begin{pmatrix} 0& 0 \\ 1 & 0 \end{pmatrix}) \cong (R \otimes _k K)_{\P}$. 
If  $x \in \P$, then  $x-t^2$  is a unit in $(R \otimes _k K)_{\P}$, hence by the same reason we have  $P_{\P} \cong (R \otimes _k K)_{\P}$. 
Therefore  $P$  is projective of rank one.

The point is to prove that $P$ is non-free. 
Assume that  $P$  is a free $R \otimes _kK$-module, then $P \cong R \otimes _k K$, since it is of rank one. 
Let  $Q$  be an $R \otimes _k V$-module defined by the same matrix  $\Phi$:
$$
\begin{CD}
(R \otimes _k V)^2 @>{\Phi}>> (R \otimes _k V)^2 @>>> Q @>>> 0. \\
\end{CD}
$$
Since the pair of matrices  $(\Phi, \Psi)$  gives a matrix-factorization of the polynomial  $x^3-y^2$  over the regular ring  $k[[x, y]] \otimes _k V$, it is easy to see that there is an exact sequence 
$$
\begin{CD}
\ldots @>{\Psi}>> (R \otimes _k V)^2 @>{\Phi}>> (R \otimes _k V)^2  @>{\Psi}>> (R \otimes _k V)^2 @>{\Phi}>> (R \otimes _k V)^2 @>{\Psi}>> \ldots 
\end{CD}
$$
In particular  $Q$  is a Cohen-Macaulay module over  $R \otimes _k V$, hence it is flat over  $V$.  
On the other hand, clearly we have  $Q_t \cong  P \cong R \otimes _k K$  and  $Q /tQ \cong \coker (\varphi) = \m$. 
Hence we must have that  $R$  degenerates to  $\m$. 
This contradicts to the following proposition. 
\end{ex}


\begin{prop}
Let  $R$  be an integral domain and $N$ a finitely generated torsion-free $R$-module. 
If $R$  degenerates to  $N$, then  $R$  is isomorphic to  $N$. 
\end{prop}

\begin{pf}
It follows from Theorem \ref{deg} that there is an exact sequence
$$
\begin{CD}
0 @>>> Z @>{\phipsi}>>  R \oplus Z @>{\fontsize{6pt}{6pt}(\alpha \  \beta)}>> N \to 0, 
\end{CD}
$$
where  $\psi$ is a nilpotent endomorphism of  $Z$. 
We take such a short exact sequence so that  $\rank \ Z$ is minimal. 
(Recall that $\rank \ Z$  is the dimension of $Z \otimes _R K_R$ over  $K_R$, where  $K_R$  is the quotient field of  $R$.) 
If  $Z =0$ then clearly $R \cong N$. 
Therefore we assume that  $Z \not= 0$, and we shall show a contradiction. 

\vspace{6pt} 
(1) If  $Z \not= 0$ then $\alpha = 0$. 
\vspace{6pt} 

In fact, if  $\alpha \not= 0$, then the mapping $\alpha : R \to N$  is an injection, since  $N$  is torsion-free over  $R$. 
In such a case, suppose  $\psi (z)=0$  for  $z \in  Z$. 
Then, since  $\alpha (\phi (z)) + \beta (\psi (z)) = 0$, we have  $\alpha (\phi (z)) = 0$ hence  $\phi (z)=0$. 
Thus we have  $\phipsi (z) =0$. 
Since  $\phipsi$ is an injection, we have  $z =0$. 
This is true for any $z \in \ker (\psi)$, hence $\psi$  must be an injection. 
However, since  $\psi$  is a nilpotent endomorphism of  $Z$, we must have   $Z =0$.     

\vspace{6pt} 
(2) The restriction $\phi |_{\ker (\psi)} : \ker (\psi) \to R$ of  $\phi$  is a surjection. 
\vspace{6pt}

Let  $r \in R$  be any element. 
Since we have shown that  $\alpha =0$, the element  $\binom {r}{0}$  of  $R \oplus Z$  belongs to $\ker (\alpha, \beta)$. 
Thus there is an element  $z \in Z$  such that  $\phipsi (z) = \binom{r}{0}$, hence 
$z \in \ker (\psi)$  and  $\phi (z) =r$. 

\vspace{6pt} 
By claim (2) we have a homomorphism  $\lambda : R \to \ker (\psi)$ such that  $\phi \cdot \lambda = 1_R$. 
Since  $\psi \cdot \lambda = 0$, we have a commutative diagram 
$$
\begin{CD}
R   @= R  \\
 @V{\lambda}VV   @V{\binom{1}{0}}VV   \\
Z  @>{\phipsi}>>  R \oplus Z. \\
\end{CD}
$$
It then follows that there is a commutative diagram with exact rows and columns: 
$$
\begin{CD}
@. 0 @. 0 @. \\
@. @VVV @VVV @. \\
@. R   @= R  \\
@.  @V{\lambda}VV   @V{\binom{1}{0}}VV   \\
0 @>>> Z  @>{\phipsi}>>  R \oplus Z @>>> N @>>> 0 \\
@.  @V{\pi}VV   @V{\fontsize{6pt}{6pt} (0 \ 1)}VV  @| @. \\
0 @>>> Z' @>{f}>>  Z @>>> N @>>> 0 \\
@. @VVV @VVV @. \\
@. 0 @. 0,  @. \\
\end{CD}
$$
where $Z' =  Z/\lambda (R)$. 
Since $\lambda$ is a splitting monomorphism, we have  $Z \cong R \oplus Z'$ and the third row can be written as follows: 
$$
\begin{CD}
0 @>>> Z' @>{\fontsize{6pt}{6pt}\binom{\phi '}{\psi '}}>>  R \oplus Z' @>>> N @>>> 0, \\
\end{CD}
$$  
where  $\psi ' = \pi \cdot f$. 
On the other hand it follows from the commutative diagram above that  $f \cdot \pi = \psi$. 
Therefore  $(\psi ') ^{n+1} = (\pi\cdot f)\cdot(\pi\cdot f)\cdots (\pi\cdot f) = \pi \cdot \psi ^n \cdot f$ for any  $n \in \N$. 
Thus  $\psi '$ is nilpotent as well as  $\psi$. 
Since  $\rank \  Z' < \rank  Z$, 
this contradicts the choice of  $Z$. 
\qed
\end{pf}

\section{Definition and properties of stable degeneration}

Let  $A$  be a commutative Gorenstein ring which is not necessarily local. 
We say that a finitely generated $A$-module  $M$ is  Cohen-Macaulay if  $\Ext _A ^i (M, A) = 0$  for all $i>0$. 
We consider the category of all Cohen-Macaulay modules over  $A$  with all $A$-module homomorphisms: 
$$
\CM (A) := \{ M \in \mod (A) \ | \ \text{$M$  is a Cohen-Macaulay module over $A$} \}, 
$$
where  $\mod (A)$ denotes the category of all finitely generated $A$-modules. 
We can then consider the stable category of  $\CM (A)$, which we denote by   
$\sCM (A)$. 
Recall that the objects of $\sCM (A)$ is  Cohen-Macaulay modules over $A$, and  the morphisms of $\sCM (A)$ are elements of 
$
\sHom _A (M, N) := \Hom _A (M, N)/P(M, N)
$ 
for  $M, N \in \sCM (A)$, where  $P(M, N)$ denotes the set of morphisms from $M$ to $N$ factoring through projective $A$-modules. 
For a  Cohen-Macaulay module  $M$ we denote it by $\u{M}$ to indicate that it is an object of  $\sCM (A)$. 
Note that  $\u{M} \cong \u{N}$  in  $\sCM (A)$  if and only if   there are projective $A$-modules  $P_1$  and  $P_2$  such that  $M \oplus P_1 \cong  N \oplus P_2$ in  $\CM (A)$.

Under such circumstances it is known that  $\sCM (A)$  has a structure of triangulated category. 
In fact, if  $L \in \CM(A)$ then we can embed  $L$  into a projective $A$-module $P$  such that the quotient  $P/L$, which we denote by  $\Omega ^{-1}L$, is Cohen-Macaulay as well. 
We define the shift functor in  $\sCM (A)$  by  $\u{L}[1] = \u{\Omega ^{-1}L}$. 
If there is an exact sequence  $0 \to L \to M \to N \to 0$  in $\CM (A)$, 
then we have the following commutative diagram by taking the pushout:
$$
\begin{CD}
0 @>>> L @>>> M @>>> N @>>> 0 \\
@. @|  @VVV  @VVV @. \\
0 @>>> L @>>> P @>>> \Omega^{-1}L @>>> 0 \\
\end{CD}
$$
We define the triangles in $\sCM(A)$  are the sequences 
$$
\begin{CD}
\u{L} @>>> \u{M} @>>> \u{N} @>>> \u{L}[1] 
\end{CD}
$$
obtained in such a way. 

Let $x \in A$  be  a non-zero divisor on  $A$. 
Note that  $x$  is a non-zero divisor on every  Cohen-Macaulay module over $A$. 
Thus the functor  $- \otimes _A A/xA$  sends a  Cohen-Macaulay module over $A$  to that over  $A/xA$. 
Therefore it yields a functor  $\CM (A) \to \CM (A/xA)$. 
 Since this functor maps projective $A$-modules to projective $A/xA$-modules, 
 it induces the functor  $\mathcal{R} : \sCM (A) \to \sCM (A/xA)$. 
It is easy to verify that  $\mathcal{R}$  is a triangle functor. 

Now let  $S \subset A$ be a multiplicative subset of  $A$. 
Then, by a similar reason to the above, we have a triangle functor  $\mathcal{L} : \sCM (A) \to \sCM (S^{-1}A)$  which maps  $\u{M}$  to  $\u{S^{-1}M}$.

\medskip

As before, let $(R, \m, k)$ be a Gorenstein local ring that is a $k$-algebra and let  $V = k[t]_{(t)}$  and  $K = k(t)$. 
Since  $R \otimes _k V$ and $R \otimes _k K$  are Gorenstein rings, we can apply the observation above. 
Actually, $t \in R\otimes _kV$  is a non-zero divisor on  $R\otimes _kV$ and 
there are isomorphisms of $k$-algebras; 
$(R\otimes _kV)/t( R\otimes _kV) \cong R$  and  $( R\otimes _kV)_t \cong  R\otimes _kK$. 
Thus there are triangle functors  $\mathcal{L} : \sCM (R\otimes _k V) \to  \sCM (R \otimes _k K)$ defined by the localization by $t$, and 
$\mathcal{R} : \sCM (R \otimes _k V) \to \sCM (R)$  defined by taking  $- \otimes _{R\otimes _kV}(R\otimes _kV)/t( R\otimes _kV) = - \otimes _V V/tV$. 
Now we define the stable degeneration of  Cohen-Macaulay modules.

\begin{defn}\label{def stable deg}
Let  $\u{M},\ \u{N} \in \sCM (R)$. 
We say that  $\u{M}$ {\bf stably degenerates to} $\u{N}$ 
if there is a  Cohen-Macaulay module  $\u{Q} \in \sCM (R \otimes _k V)$  
such that $\mathcal{L}( \u{Q} ) \cong \u{ M \otimes _k K}$  in  $\sCM (R\otimes _k K)$ and $\mathcal{R} (\u{Q}) \cong \u{N}$  in $\sCM (R)$. 
\end{defn}

The following are easily verified. 
\begin{lemma}
Let  $M, N \in \CM (R)$. 
If  $M$  degenerates to $N$, then  $\u{M}$  stably degenerates to $\u{N}$.   
\end{lemma} 

\begin{pf}
Suppose  $M$  degenerates to $N$. 
Then there is  $Q \in \mod (R \otimes _k V)$ that is $V$-flat and satisfies 
$Q/tQ \cong N$  and  $Q_t \cong M \otimes _k K$. 
Note in this case that $Q$  is a Cohen-Macaulay $R\otimes _kV$-module, i.e. $Q \in \CM (R \otimes _kV)$. 
In fact, if a prime ideal $\P$ of  $R \otimes _kV$  contains $t$, 
then, since  $Q_{\P}/tQ_{\P}=(Q/tQ)_{\P}$ is Cohen-Macaulay and since  $t$  is a non-zero divisor on  $Q$, we see that  $Q_{\P}$  is a Cohen-Macaulay $R_{\P}$-module. 
On the other hand, if  $t \not\in \P$, then $Q_{\P}$  is a Cohen-Macaulay $R_{\P}$-module as well, since  it is a localization of the Cohen-Macaulay $R \otimes _k K$-module  $Q_t$. 
Thus we have  $Q \in \CM (R \otimes _k V)$. 
Then it is clear that  $\mathcal{R}(\u{Q}) \cong \u{N}$  and  $\mathcal{L}(\u{Q}) \cong \u{M \otimes _k K}$. 
\qed
\end{pf}

\begin{prop}\label{triangle degeneration}
Suppose that there is a triangle  
$$
\begin{CD}
\u{L} @>{\alpha}>> \u{M} @>{\beta}>> \u{N} @>{\gamma}>> \u{L}[1], 
\end{CD}
$$
in $\sCM (R)$. 
Then  $\u{M}$  stably degenerates to  $\u{L} \oplus \u{N}$. 
\end{prop}

\begin{pf}
Note that $\gamma$  is an element of  $\sHom _{R}(N, \Omega^{-1}L)$, which is naturally a submodule of $\sHom _{R}(N, \ \Omega^{-1}L)\otimes _k V \cong \sHom _{R\otimes _k V}(N\otimes _k V, \ \Omega^{-1}L\otimes _k V)$. 
Regarding  $\gamma$  as an element of  $\sHom _{R\otimes _k V}(N\otimes _k V, \ \Omega^{-1}L\otimes _k V)$, we have a morphism 
$\gamma \otimes t : \u{N \otimes _k V} \to \u{L\otimes _k V}[1]$  in $\sCM (R\otimes _k V)$. 
(The morphism $\gamma \otimes t$  is a composition of  $\gamma :  \u{N \otimes _k V} \to \u{L\otimes _k V}[1]$  with the multiplication map by $t$.)
Now embed $\gamma \otimes t$ into a triangle in $\sCM (R\otimes _k V)$, and we get $\u{Q}\in \sCM (R\otimes _k V)$  with the triangle 
$$
\begin{CD}
\u{L\otimes _k V} @>>> \u{Q}  @>>>  \u{N \otimes _k V} @>{\gamma \otimes t}>>  \u{L\otimes _k V}[1].
\end{CD}
$$
Notice that  $\mathcal{R} (\gamma \otimes t) = 0$, since the multiplication map by $t$  is zero in $\sCM (R)$.  
Therefore applying the triangle functor $\mathcal{R}$  to the triangle above and noting that  $\mathcal{R}(\u{X\otimes _k V}) = \u{X}$  for any  $\u{X} \in \sCM (R)$, we have an isomorphism $\mathcal{R}(\u{Q}) \cong \u{L} \oplus \u{N}$. 

On the other hand, there is an isomorphism of triangles in  $\sCM(R\otimes _k K)$; 
$$
\begin{CD}
\u{L\otimes _k K} @>{\alpha \otimes 1}>> \u{M\otimes _kK}  @>{\beta \otimes 1}>>  \u{N \otimes _k K} @>{\gamma \otimes 1}>>  \u{L\otimes _k K}[1] \\
@V{t}V{\cong}V @VVV @| @V{t}V{\cong}V \\
\u{L\otimes _k K} @>>> \u{Q_t}  @>>>  \u{N \otimes _k K} @>{\gamma \otimes t}>>  \u{L\otimes _k K}[1].
\end{CD}
$$
Hence we have that $\mathcal{L}(\u{Q}) = \u{Q_t} \cong \u{M \otimes _k K}$. 
As a consequence,  $\u{M}$  stably degenerates to  $\u{L}\oplus \u{N}$.
\qed
\end{pf}

\begin{prop}\label{shift and dual}
Let  $\u{M}, \ \u{N} \in \u{\rm CM}(R)$  and suppose that  $\u{M}$  stably degenerates to $\u{N}$. Then the following hold. 
\begin{itemize}
\item[$(1)$]
$\u{M}[1]$ (resp. $\u{M}[-1]$) stably degenerates to  $\u{N}[1]$ (resp. $\u{N}[-1]$). 
\item[$(2)$]
$\u{M^*}$ stably degenerates to  $\u{N^*}$, where  $M^*$ denotes the $R$-dual $\Hom _R (M, R)$. 
\end{itemize}
\end{prop}

\begin{pf}
By the assumption, there is $\u{Q} \in \sCM(R \otimes _k V)$  
such that $\mathcal{L}( \u{Q} ) \cong \u{ M \otimes _k K}$  in  $\sCM (R\otimes _k K)$ and $\mathcal{R} (\u{Q}) \cong \u{N}$  in $\sCM (R)$. 
 
To prove (1), consider  $\u{Q}[1]$  in  $\sCM(R\otimes _kV)$. Then, since  $\mathcal{L}$  and  $\mathcal{R}$  are triangle functors, we have  $\mathcal{L}( \u{Q}[1] ) \cong \u{ M \otimes _k K}[1]\cong \u{ M [1] \otimes _k K}$  in  $\sCM (R\otimes _k K)$ and $\mathcal{R} (\u{Q}[1]) \cong \u{N}[1]$  in $\sCM (R)$. 
This shows that  $\u{M}[1]$  stably degenerates to  $\u{N}[1]$. 
 
To prove (2) we set  $\tilde{Q} = \Hom _{R\otimes _kV} (Q, \ R \otimes _kV)$, and consider  $\u{\tilde{Q}} \in \sCM (R \otimes _kV)$. 
Then it is easy to see that there are isomorphisms; 
$\tilde{Q}_t \cong  \Hom _{R\otimes _kK} (Q_t, \ R \otimes _kK)$
and 
$\tilde{Q}/t \tilde{Q} \cong \Hom _R (Q/tQ, R)$. 
Thus  $\mathcal{L}(\u{\tilde{Q}}) \cong \u{\Hom _{R\otimes _kK}(M \otimes _k K, R\otimes _k K)} \cong \u{\Hom _R (M, R) \otimes _kK}$  and 
$\mathcal{R}(\u{\tilde{Q}})\cong \u{\Hom _R(N, R)}$. 
Therefore  $\u{M^*}$ stably degenerates to  $\u{N^*}$. 
\qed
\end{pf}

\begin{prop}\label{strange degeration}
Let  $\u{M}, \ \u{N}, \ \u{X} \in \u{\rm CM}(R)$. 
If  $\u{M}\oplus \u{X}$  stably degenerates to  $\u{N}$, 
then  $\u{M}$  stably degenerates to  $\u{N}\oplus \u{X}[1]$. 
\end{prop}

\begin{pf}
Take $\u{Q} \in \sCM(R \otimes _k V)$  satisfying that 
$\mathcal{L}( \u{Q} ) \cong \u{ (M \oplus X) \otimes _k K}$  in  $\sCM (R\otimes _k K)$ and $\mathcal{R} (\u{Q}) \cong \u{N}$  in $\sCM (R)$. 
 
First of all we note that there is a natural isomorphism 
$$
\sHom _{R \otimes _k V} ((M \oplus X) \otimes _k V,\ Q)_t \cong 
\sHom _{R \otimes _k K} ((M \oplus X) \otimes _k K,\ Q_t).
$$
Thus there is a morphism  $\alpha : \u{(M \oplus X) \otimes _k V} \to \u{Q}$  in  $\sCM (R \otimes _kV)$ such that  $\mathcal{L} (\alpha) = \alpha_t$  is an isomorphism in  $\sCM (R \otimes _kK)$.
Here, replacing  $\alpha$  with $t \alpha$  if necessary, we may assume that  $\mathcal{R}(\alpha) =0$.  
Now let  $\beta : \u{X \otimes _k V} \to \u{(M \otimes _k V) \oplus (X \otimes _k V)}\cong \u{(M \oplus X) \otimes _k V}$ be a natural splitting monomorphism and we set  $\gamma = \alpha \cdot \beta$. 
Then we can embed the morphism $\gamma$ into a triangle in  $\sCM (R \otimes _kV)$: 
$$
\begin{CD}
\u{X \otimes _k V} @>{\gamma}>> \u{Q} @>>> \u{Q'} @>>> \u{X \otimes _k V}[1]. 
\end{CD}
$$
Since  $\mathcal{R}(\gamma) = \mathcal{R}(\alpha)\cdot \mathcal{R}(\beta) =0$, we have that 
$$
\mathcal{R}(\u{Q'}) \cong \mathcal{R}(\u{Q}) \oplus \mathcal{R}(\u{X\otimes _k V}[1]) \cong \u{N} \oplus \u{X}[1]. 
$$ 
On the other hand, since there is a triangle in $\sCM (R \otimes _kK)$; 
$$   
\begin{CD}
\u{X \otimes _k K} @>{\mathcal{L}(\beta)}>> \u{(M \oplus X) \otimes _k K} @>>> \u{M \otimes _k K} @>>> \u{X \otimes _k K}[1], 
\end{CD}
$$
noting that  $\mathcal{L}(\alpha): \u{(M \oplus X) \otimes _k K} \to \u{Q_t}$  is an isomorphism, we have that $\mathcal{L}(\u{Q'}) \cong \u{M\otimes _kK}$. 
This shows that   $\u{M}$  stably degenerates to  $\u{N}\oplus \u{X}[1]$. 
\qed
\end{pf}

\begin{rem}
The zero object in $\sCM (R)$  can stably degenerate to a non-zero object. 
For example,  there is a triangle 
$$
\begin{CD}
\u{X} @>>> \u{0} @>>> \u{X}[1] @>{1}>> \u{X}[1]. 
\end{CD}
$$
for any $\u{X} \in \sCM (R)$.
Hence $\u{0}$  stably degenerates to $\u{X} \oplus \u{X}[1]$ by Proposition \ref{triangle degeneration}.   
\end{rem}


\section{Conditions for stable degeneration}


Let  $(R, \m, k)$ be a Gorenstein local $k$-algebra as before. 
The main purpose of this section is to prove the following theorem.

\begin{thm}\label{implication}
Consider the following three conditions for  Cohen-Macaulay $R$-modules $M$  and  $N$ : 

\begin{itemize}
\item[$(1)$]
$R ^m \oplus M$  degenerates to  $R^n \oplus N$  for some  $m, n \in {\Bbb N}$. 

\item[$(2)$]
There is a triangle  
$\u{Z} \overset{\binom{\u{\phi}}{\u{\psi}}}\to \u{M} \oplus \u{Z} \to \u{N} \to \u{Z}[1]$  in  $\u{\rm CM}(R)$,  where $\u{\psi}$  is  a nilpotent element of  $\u{\rm End} _R (Z)$. 

\item[$(3)$]
$\u{M}$ stably degenerates to  $\u{N}$. 
\end{itemize}

\noindent
Then, in general, the implications $(1) \Rightarrow (2) \Rightarrow (3)$ hold. 

\noindent
Furthermore, if $R$  is artinian, then the conditions  $(1), (2)$  and  $(3)$  are all equivalent. 
\end{thm}

In the next section we shall prove that  the implication  $(3) \Rightarrow (2)$  holds if  $R$  is a Gorenstein complete local ring with only an isolated singularity (Theorem (\ref{isolated singularity})).

\begin{pf}
$(1) \Rightarrow (2)$: 
Suppose that $R ^m \oplus M$  degenerates to  $R^n \oplus N$. 
Then by Theorem \ref{deg}, we have a short exact sequence of finitely generated left $R$-modules
$$
0 \to Z \overset{\phipsi}\longrightarrow (R^m\oplus M) \oplus Z \to (R^n \oplus N) \to 0,
$$
where $\psi$ is nilpotent. 
In such a case $Z$ ia a  Cohen-Macaulay module as well. 
See \cite[Remark 4.3.1]{Y3} or \cite[Proof of Theorem 3.2]{Y4}. 
Then converting this short exact sequence into a triangle in $\sCM (R)$, we have $\u{Z} \overset{\binom{\u{\phi}}{\u{\psi}}}\to \u{M} \oplus \u{Z} \to \u{N} \to \u{Z}[1]$, where, since  $\psi \in \End _R (Z)$  is nilpotent, 
$\u{\psi} \in \sEnd _R (Z)$  is nilpotent as well.

\medskip
$(2) \Rightarrow (3)$: 
Through the natural injective homomorphism 
$$
\begin{array}{lrll}
&\sHom _R (Z, M) &\hookrightarrow &\sHom _{R\otimes _k V} (Z \otimes _k V, M \otimes _k V), \vspace{6pt} \\
&\u{\phi}  &\mapsto  &\u{\phi \otimes 1_V}  
\end{array}
$$
we regard $\u{\phi}$ as a morphism in  $\sCM (R \otimes _k V)$. 
Likewise $\u{\psi}$  is regarded as a morphism in $\sCM (R \otimes _k V)$  which is nilpotent as well. 
Note that  there is a triangle in $\sCM (R \otimes _kV)$; 
$$
\u{Z\otimes _k V} \overset{\binom{\u{\phi}}{\u{\psi}}}\longrightarrow  \u{M\otimes _k V} \oplus \u{Z\otimes _k V} \longrightarrow  \u{N\otimes _k V} \longrightarrow  \u{Z\otimes _k V}[1].
$$ 
Now consider a morphism  $t + \u{\psi} : \u{Z\otimes _k V} \to \u{Z\otimes _k V}$, and we have a triangle of the form; 
$$
\u{Z\otimes _k V} \overset{\binom{\u{\phi}}{t+\u{\psi}}}\longrightarrow  \u{M\otimes _k V} \oplus \u{Z\otimes _k V} \longrightarrow  \u{Q} \longrightarrow  \u{Z\otimes _k V}[1], 
$$
for some  $\u{Q} \in \sCM (R\otimes _kV)$.
Note that  $\mathcal{L}(t + \u{\psi})$  is an isomorphism in  $\sCM (R \otimes _k K)$, since  $t \in R\otimes _k K$  is a unit and  $\mathcal{L}(\u{\psi})$  is a nilpotent morphism. 
Thus, applying the functor $\mathcal{L}$  to the triangle above, we have 
that  $\mathcal{L}(\u{Q}) \cong \mathcal{L}(\u{M \otimes _kV}) = \u{M \otimes _kK}$. 
On the other hand, since  $\mathcal{R}(t + \u{\psi})= \u{\psi}$, we see that 
$\mathcal{R}(\u{Q}) \cong N$. 
Thus  $\u{M}$  stably degenerates to  $\u{N}$.

\medskip
$(3) \Rightarrow (1)$: 
In this proof we assume that  $\dim R =0$. 
Since $\u{M}$ stably degenerates to  $\u{N}$, there is $\u{Q} \in \sCM(R \otimes _k V)$  such that $\mathcal{L}( \u{Q} ) \cong \u{ M \otimes _k K}$  in  $\sCM (R\otimes _k K)$ and $\mathcal{R} (\u{Q}) \cong \u{N}$  in $\sCM (R)$.
By definition, we have isomorphisms 
$Q_t \oplus P_1 \cong (M \otimes _k K) \oplus P_2$ in $\CM (R\otimes _kK)$ 
for some projective  $R\otimes _k K$-modules $P_1, P_2$,  and 
$Q/tQ \oplus R^a \cong N \oplus R^b$ in $\CM (R)$  for some $a, b \in \N$. 
As we have remarked in Remark \ref{remark on RotimesK}, $R \otimes _k K$  is a local ring, hence  $P_1$  and  $P_2$  are free. 
Thus 
$Q_t \oplus (R\otimes _kK)^c  \cong (M \otimes _k K) \oplus (R\otimes _kK)^d$ 
for some $c, d \in \N$. 

Now setting  $\tilde{Q} = Q \oplus (R\otimes _kV)^{a+c}$, we have isomorphisms
$$
\tilde{Q}_t \cong (M \oplus R^{a+d}) \otimes _k K, \quad 
\tilde{Q}/t\tilde{Q} \cong N \oplus R^{b+c}.
$$
Noting that  $Q$, hence $\tilde{Q}$, is $V$-flat, since it is a  Cohen-Macaulay module over  $R\otimes _kV$, 
we conclude that $M \oplus R^{a+d}$ degenerates to $N \oplus R^{b+c}$. 
\qed
\end{pf}

In Theorem \ref{implication} the implication  $(2) \Rightarrow (1)$  does not holds even in the case when  $\dim  R =1$. 

\begin{ex}\label{ex 2->1}
As in Example \ref{nonfree-projective}, let  $R = k[[ x, y]]/(x^3-y^2)$  and let $\m = (x, y)R$. 
The ring $R$ can be identified with the subring  $k[[s^2, s^3]]$  of the formal power series ring  $S=k[[s]]$ by mapping  $x, y$ to  $s^2, s^3$ respectively. 
Note in this case that  $\End _R (\m) = S$. 
Actually  $\m$  is the conductor of the ring extension  $R \subset S$ and hence  $\m$  is identical with an ideal $s^2S$  of  $S$.

In this case, there is an exact sequence 
$$
\begin{CD}
0 @>>> \m @>{\binom{j}{s}}>> R \oplus \m @>{(x \ -s)}>> \m @>>> 0,  \\
\end{CD}
$$
where  $j$  is a natural inclusion  $\m \subset R$. 
This is in fact a unique AR sequence in $\CM (R)$. (See \cite[Proposition 5.11]{Y1}.) 
Note that  $\sEnd _R (\m) = S/s^2S$  and there is a triangle 
$$
\begin{CD}
\u{\m} @>{\u{s}}>> \u{\m} @>\u{-s}>> \u{\m} @>>> \u{\m}[1],  \\
\end{CD}
$$
in $\sCM (R)$. 
Since  $\u{s} \in  \sEnd _R (\m)=S/s^2S$  is nilpotent, the condition $(2)$ in Theorem \ref{implication} holds for  $M = 0$  and  $N = \m$. 
Hence  $\u{0}$  stably degenerates to  $\u{\m}$  in this case. 
However we can prove the following proposition, and hence the condition $(1)$  in Theorem \ref{implication} does not hold. 
\end{ex}

\begin{prop}
Let $R = k[[s^2, s^3]] \subset S = k[[s]]$  and let $\m = (s^2, s^3)R$ as above. Then $R ^m$  never degenerates to  $R^n \oplus \m$  for  $m, n \in {\Bbb N}$.  
\end{prop}

\begin{pf}
Suppose that $R^m$  degenerates to  $R^n \oplus \m$  for some  $m, n \in \N$, and we shall seek a contradiction. 
Comparing the ranks, we have  $m = n+1$  by Remark \ref{if degerates}(1). 
Then it follows from Theorem \ref{deg} that there is an exact sequence 
$$
0 \to Z \overset{\phipsi}\longrightarrow R^{n+1} \oplus Z \to R^{n} \oplus \m \to 0,
$$
such that the endomorphism  $\psi$  of  $Z$ is nilpotent. 
Note that  $Z$ is also a Cohen-Macaulay module over  $R$. 
Since  $R$  and  $\m$  are unique indecomposable  Cohen-Macaulay $R$-modules over  $R$,  $Z$ can be described as  $R ^a \oplus \m^b$  for some $a, b \in \N$. 
Therefore the above sequence is described as 
\def\mat{\begin{pmatrix}
\phi_1 & \phi_2 \\
\psi_{11} & \psi_{12} \\
\psi_{21} & \psi_{22} \\
\end{pmatrix}}
$$
\begin{CD}
0 @>>> R^a \oplus \m ^b  @>{\mat}>> R^{n+1} \oplus R^a \oplus \m ^b  @>>> R^{n} \oplus \m @>>> 0,  \\
\end{CD}
$$
where  $\phi = ( \phi_1 \ \phi _2)$  and  
$\psi = \begin{pmatrix}
\psi_{11} & \psi_{12} \\
\psi_{21} & \psi_{22} \\
\end{pmatrix}$.
Note that  $S = \End _R (\m) = \Hom _R (\m, R)$  and  $\Hom _R (R, \m) = \m$. 
Thus, 
$$
\begin{pmatrix}
\psi_{11} & \psi_{12} \\
\psi_{21} & \psi_{22} \\
\end{pmatrix}
\in 
\begin{pmatrix}
R ^{a \times a} & S^{a \times b} \\
\m ^{b \times a} & S^{b \times b} \\
\end{pmatrix}
$$
As the first step of the proof we claim that 
\begin{equation}\label{step1}
\det (\psi _{22}) \equiv c s \ \ \ (\mod \ s^2S), 
\end{equation}
for some $c \in k\backslash \{0\}$.

In fact, converting the short exact sequence above into the stable category, we have the triangle in  $\sCM (R)$; 
$$
\begin{CD}
\u{\m}^b @>{\u{\psi}}>> \u{\m}^b @>>> \u{\m} @>>> \u{\m}^b[1],  \\
\end{CD}
$$
where we should note that  $\u{\psi} = \u{\psi_{22}}$. 
Note that $\sEnd _R (\m) = S/s^2S$  and  $\u{\m} [1] \cong \u{\m}$ in $\sCM (R)$. 
It is easy to see that there are triangles in $\sCM (R)$; 
$$
\begin{cases} 
\ \ \ \u{\m} \overset{1}\longrightarrow \u{\m} \longrightarrow  \u{0} \longrightarrow \u{\m}, \\
\ \ \ \u{\m} \overset{\u{s}}\longrightarrow  \u{\m} \longrightarrow  \u{\m} \longrightarrow  \u{\m}, \\
\ \ \ \u{\m} \overset{0}\longrightarrow  \u{\m} \longrightarrow  \u{\m}^{\oplus 2} \longrightarrow \u{\m}. \\
\end{cases}
$$
Since  $\u{\psi _{22}}$  is an element  $(S/s^2S)^{b \times b}$,  noticing that the cone of  $\u{\psi _{22}}$  is $\u{\m}$, we can make  $\u{\psi _{22}}$  into the following form after elementary transformations of matrices over  $S/s^2S$:
$$
\begin{pmatrix}
1 & &  &  \\
  &\ddots &  &  \\
  & & 1&  \\
  & &  &s \\
\end{pmatrix}
$$
In particular, $\det (\u{\psi_{22}})$  is equal to  $s$  up to a unit in  $S/s^2S$. 
Since the natural projection  $S \to S/s^2S$  sends  $\det (\psi_{22})$ to  $\det (\u{\psi_{22}})$, this shows the equality (\ref{step1}). 

Now let us denote by  $\chi _{\psi_{22}}(T)$ the characteristic polynomial of the matrix  $\psi _{22} \in S^{b \times b}$. 
We have from (\ref{step1}) that 
\begin{equation}\label{step2}
\chi _{\psi_{22}}(T) = \det (TE - \psi_{22}) \equiv T^b + \cdots + (-1)^b cs \ \ \ (\mod \ s^2S).
\end{equation}
 
Now we consider the matrix 
$\psi = 
\begin{pmatrix}
\psi_{11} & \psi_{12} \\
\psi_{21} & \psi_{22} \\
\end{pmatrix}
$ as an element of 
$
\begin{pmatrix}
S ^{a \times a} & S^{a \times b} \\
S ^{b \times a} & S^{b \times b} \\
\end{pmatrix}, 
$
which is a nilpotent matrix as well. 
Therefore, as for the characteristic polynomial  $\chi _{\psi} (T)$ of $\psi$  we have the equality 
\begin{equation}\label{step3}
\chi _{\psi} (T) = T^{a+b}.
\end{equation}
On the other hand, since every entry of  $\psi _{21}$ is in $\m = s^2S$, we have 
$$
\psi \equiv  
\begin{pmatrix}
\psi_{11} & \psi_{12} \\
0 & \psi_{22} \\
\end{pmatrix}
\ \ \ (\mod \ s^2 S), 
$$
and thus 
\begin{equation}\label{step4}
\chi _{\psi}(T) \equiv  \chi _{\psi_{11}}(T) \cdot \chi _{\psi_{22}}(T) \ \ \ (\mod \ s^2S).
\end{equation}
Since every entry of  $\psi _{11}$  is in  $R = k[[s^2, s^3]]$, 
we see that  $\psi_{11} \mod \ s^2S$  is just a matrix with entries in $k$. 
Thus  $\chi _{\psi_{11}} (T) \ \mod \ s^2S \in k[T]$. 
Combining this observation with (\ref{step2}),(\ref{step3}) and (\ref{step4}), we have the equality of elements in  $(S/s^2S)[T]$;  
$$
T^{a+b} = (T^a + d_1T^{a-1} + \cdots + d_a)(T^b + \cdots + (-1)^b cs), 
$$ 
with  $d_i \in k$. 
Setting  $d_0=1$  and  $\ell = \max \{\ i \ | \ 0 \leq i \leq a, \  d_i \not= 0 \}$, we see that the nontrivial term  $(-1)^b d_{\ell}csT^{a-\ell}$ appears in the right-hand side.  
This is a contradiction. 
\qed
\end{pf}

\section{The case of isolated singularity}

In this section we shall prove the equivalence  $(2) \Leftrightarrow (3)$  for the conditions in Theorem \ref{implication} if the Gorenstein local ring is an isolated singularity.
The goal of this section is to prove the following theorem.

\begin{thm}\label{isolated singularity}
Let  $(R, \m, k)$  be a Gorenstein complete local ring that is a $k$-algebra, 
and let  $\u{M}, \u{N} \in \sCM (R)$. 
Assume that $R$ has only an isolated singularity, and $k$  is an infinite field. If $\u{M}$ stably degenerates to  $\u{N}$, 
then there is a triangle in  $\sCM (R)$;  
$$
\begin{CD}
\u{Z} @>{\binom{\u{\phi}}{\u{\psi}}}>>  \u{M} \oplus \u{Z} @>>> \u{N} @>>>  \u{Z}[1], 
\end{CD}
$$
where $\u{\psi}$  is  a nilpotent element of  $\u{\rm End} _R (Z)$. 
\end{thm}

To prove this theorem we need some auxiliary lemmas. 
For the proof of the following lemma the reader should refer to \cite[Swan's Lemma 5.1]{L}.


\begin{lemma}[Swan's Lemma]\label{Swan}
Let  $R$  be a noetherian ring and  $t$  a variable. 
Assume that an $R[t]$-module  $L$  is a submodule of $W\otimes _R R[t]$  with  $W$ being a finitely generated  $R$-module. 
Then, there is an exact sequence of $R[t]$-modules; 
$$
\begin{CD}
0 @>>> X \otimes _R R[t] @>>> Y \otimes _R R[t] @>>> L @>>> 0, 
\end{CD}
$$
where  $X$  and $Y$  are finitely generated $R$-modules.
\end{lemma}


\begin{lemma}\label{Swan lemma}
Let  $(R, \m, k)$  be a noetherian local $k$-algebra, where  $k$  is an infinite field, and let  $V = k[t]_{(t)}$  and  $K=k(t)$  as before. 
Suppose that a finitely generated  $R \otimes _k V$-module  $P'$  satisfies the following conditions: 
\begin{itemize}
\item[$(1)$]  
$P'$  is a submodule of a finitely generated free $R\otimes _k V$-module,  
\item[$(2)$]
the localization  $P = P'_t$  by  $t$  is a projective  $R \otimes _k K$-module. 
\end{itemize}
Then there is a short exact sequence of  $R \otimes _k V$-modules; 
$$
\begin{CD}  
0 @>>> X \otimes _kV @>>> (X \otimes _k V)  \oplus  (R ^n \otimes _k V) @>>> P' @>>> 0, 
\end{CD}
$$
where   $X$  is a finitely generated  $R$-module and  $n$ is a non-negative integer. 
\end{lemma}

\begin{pf}
Recall that  $R \otimes _k V = S^{-1} R[t]$  and  $R \otimes _k K = T^{-1} R[t]$  where  $S = k[t] \backslash (t)$  and  $T = k[t] \backslash \{ 0\}$. 
Since  $P'$  is a finitely generated  $S^{-1}R[t]$-module, we find a finitely generated $R[t]$-submodule  $L$ of  $P'$  which satisfies  $S^{-1} L = P'$. 
By the assumption, we may assume $P' \subseteqq S^{-1}R[t] ^r$ for some $r \geqq 0$. 
Replacing  $L$  with  $L \cap R[t]^r$ if necessary, we may take  a submodule of a free $R[t]$-module as such an $L$. 
Thus we can apply Swan's Lemma \ref{Swan} to  $L$ and we have an exact sequence of $R[t]$-modules
\begin{equation}\label{s1}
\begin{CD}
0 @>>> X \otimes _R R[t] @>>> Y \otimes _R R[t] @>>> L @>>> 0,
\end{CD}
\end{equation}
where  $X, Y$  are finitely generated $R$-modules.

We note that there is a polynomial  $f(t) \in T$  such that  $L_{f(t)}$  is projective over $R[t]_{f(t)}$. 
In fact, taking an epimorphism from a free module to  $L$, we have an exact sequence  $0 \to L_1 \to R[t]^m \to L \to 0$. 
This sequence will split if it is localized by $T$, since  $P = T^{-1}L$  is projective. 
Therefore it splits after localizing it by some element  $f(t) \in T$.

Now we localize the exact sequence  (\ref{s1})  by  $f(t) \in T$, and then we have an isomorphism of  $R[t]_{f(t)}$-modules 
\begin{equation}\label{s2} 
Y \otimes _R R[t]_{f(t)} \cong L_{f(t)} \oplus  (X \otimes _R R[t]_{f(t)}). 
\end{equation}
Since we assume that  $k$  is an infinite field, there is an element  $c \in k$  such that  $f(c) \not=0$. 
Consider the $R$-algebra homomorphism  $\sigma : R[t]_{f(t)} \to  R$  which sends $t$  to  $c$, and we regard  $R$  as an $R[t]_{f(t)}$-algebra through $\sigma$. 
We note that   $(X \otimes _R R[t]_{f(t)})\otimes  _{R[t]_{f(t)}} R \cong X$ for any $R$-module  $X$.
Hence, taking   $- \otimes  _{R[t]_{f(t)}} R$ with  (\ref{s2}), we have  
$$
Y \cong (L_{f(t)} \otimes _{R[t]_{f(t)}} R) \oplus X. 
$$
Noting that $L_{f(t)}$  is a projective  $R[t]_{f(t)}$-module, we see that 
$L_{f(t)} \otimes _{R[t]_{f(t)}} R$  is a projective  $R$-module as well.
Hence  $L_{f(t)} \otimes _{R[t]_{f(t)}} R \cong R^n$  as $R$-modules for some $n \geq 0$, since  $R$  is local. 
Therefore we have  $Y \cong  R^n \oplus X$.
Substituting this into  (\ref{s1})   and taking the localization by $S$, we have a short exact sequence 
$$
\begin{CD}  
0 @>>> X \otimes _kV @>>> (X \otimes _k V)  \oplus  (R ^n \otimes _k V) @>>> P' @>>> 0. 
\end{CD} 
$$
\qed
\end{pf}


\begin{prop}\label{Swan triangle}
Let  $R$  be a Gorenstein local $k$-algebra, where  $k$  is an infinite field. 
Suppose we are given a  Cohen-Macaulay  $R \otimes _k V$-module  $P'$  
satisfying that the localization  $P = P'_t$  by  $t$  is a projective  $R \otimes _k K$-module. 
Then 
there is a Cohen-Macaulay  $R$-module  $X$  with a triangle in  $\sCM (R \otimes _k V)$  of the following form: 
\begin{equation}\label{short exact}
\begin{CD}  
\u{X \otimes _kV} @>>> \u{X \otimes _k V} @>>> \u{P'} @>>> \u{X \otimes _k V} [1]. 
\end{CD}
\end{equation}
\end{prop}

\begin{pf}
Let $\ell = \dim R +1 = \dim R \otimes _k V$, and we note that $P'$  is an $\ell$th syzygy module of a  Cohen-Macaulay $R\otimes _k V$-module  $P''$. 
($R\otimes _k V$  is a Gorenstein ring and every  Cohen-Macaulay module over a Gorenstein ring is a syzygy module of a  Cohen-Macaulay module.) 
In such a situation, $P''_t$  is a projective  $R\otimes _kK$-module, since so is  $P'_t$. 
Therefore we can apply Lemma \ref{Swan lemma} to  $P''$  to get a short exact sequence
$$
\begin{CD}  
0 @>>> X \otimes _kV @>>> (X \otimes _k V)  \oplus  (R ^n \otimes _k V) @>>> P'' @>>> 0,  
\end{CD}
$$
where $X$ is just a finitely generated $R$-module.
Now take the $\ell$th syzygy of this sequence.
We notice that the $\ell$th syzygy module  $\Omega _{R\otimes_k V} ^{\ell}(X \otimes _k V)$  is isomorphic,  as an object of  $\sCM (R\otimes _kV)$,  to  $\Omega _R^{\ell}(X) \otimes _k V$  which is a  Cohen-Macaulay $R\otimes _kV$-module.
In this way we obtain a triangle in  $\sCM (R \otimes _k V)$ ;
$$
\begin{CD}  
\u{\Omega _R^{\ell}X \otimes _kV} @>>> \u{\Omega _R^{\ell} X \otimes _k V} @>>> \u{\Omega _{R\otimes _k V}^{\ell} P''} @>>> \u{\Omega _R^{\ell}X \otimes _k V} [1]. 
\end{CD}
$$
Setting $P'=\Omega_{R\otimes_kK}^\ell P''$ and replacing $\Omega_R^\ell X$ with $X$,  we have a desired triangle.
\qed
\end{pf}


The following lemma is an analogue,  or one might say a higher-dimensional version,  of the Fitting lemma.

\begin{lemma}\label{Fitting lemma for isolated sing}
Let  $R$  be a Gorenstein complete local ring which has only an isolated singularity, and let  $\u{X} \in \sCM (R)$. 
Given an endomorphism  $\u{\psi} \in \sEnd _R (X)$, we have a direct decomposition  $\u{X} = \u{X_1} \oplus \u{X_2}$ and automorphisms  $\u{\alpha}, \ \u{\beta}$  of  $\u{X}$  such that 
$$
\u{\alpha} \cdot \u{\psi} \cdot \u{\beta} =  \begin{pmatrix} \u{\psi_1} & 0 \\ 0 & \u{\psi_2} \\ \end{pmatrix},
$$
where  
$\u{\psi_1} : \u{X_1} \to \u{X_1}$  is an automorphism and  
$\u{\psi_2} : \u{X_2} \to \u{X_2}$  is a nilpotent endomorphism. 
\end{lemma}

\begin{pf}
Recall that  $\CM (R)$ is a Krull-Schmidt category and so is $\sCM(R)$. 
Therefore $\u{X}$  is uniquely decomposed into a direct sum of indecomposable objects; $\u{X} \cong \u{Y_1} \oplus \cdots \oplus \u{Y_n}$. 
According to this decomposition, $\u{\psi}$  is described as  an $n \times n$-matrix  $(\u{\psi_{ij}})$, where  $\u{\psi_{ij}} \in \sHom _R(Y_j, Y_i)$. 
If there is an isomorphism  $\u{\psi_{ij}}$, then $\u{\psi}$  is arranged into the form  $\begin{pmatrix} \u{\psi_{ij}} & 0 \\ 0 & * \end{pmatrix}$, more precisely, there are automorphisms  $\u{\alpha _1}$  and $\u{\beta _1}$ of  $\u{X}$  such that  $\u{\alpha _1} \cdot \u{\psi} \cdot \u{\beta _1} = \begin{pmatrix} \u{\psi_{ij}} & 0 \\ 0 & * \end{pmatrix}$. 
Hence, by induction on   $n$, it is enough to prove that  $\u{\psi}$  is nilpotent if all $\u{\psi_{ij}}$  are non-isomorphic.     
For an integer  $N$, each entry of the matrix  $(\u{\psi_{ij}})^N$  is a composition of morphisms ;  $\u{Y_j} = \u{Z_0} \to \u{Z_1} \to \cdots \to \u{Z_N} = \u{Y_i}$, where each $\u{Z_k}$  is one of  $\u{Y_1}, \ldots , \u{Y_n}$. 
Take an integer  $\ell$  satisfying  $(\mathrm{rad} \ \sEnd _R (Y_i))^{\ell} = 0$  for $i = 1,2 , \ldots ,n$. 
Note that this is possible, since each  $\sEnd _R (Y_i)$  is an artinian local ring by the assumption.
If  $N > \ell n$, then  some  $\u{Y_k}$  appears at least $(\ell +1)$-times among  $\u{Z_0}, \u{Z_1}, \cdots , \u{Z_N}$, therefore the composition 
 $\u{Z_0} \to \u{Z_1} \to \cdots \to \u{Z_N}$  is a zero morphism.  
Thus  $(\u{\psi_{ij}})^N = 0$  if  $N > \ell n$.  
\qed
\end{pf}

\vspace{12pt}
Now we proceed to the proof of Theorem \ref{isolated singularity}.

\vspace{12pt}
[Proof of Theorem \ref{isolated singularity}]

To prove the theorem, let  $(R, \m, k)$  be a Gorenstein complete local $k$-algebra, and assume that $R$ has only an isolated singularity and that $k$  is an infinite field. 

For  $\u{M}, \u{N} \in \sCM (R)$, we assume that $\u{M}$ stably degenerates to  $\u{N}$. 
Then, by definition we have  $\u{Q} \in \sCM (R\otimes _k V)$  such that 
$\u{Q_t} \cong \u{M \otimes _k K}$  in   $\sCM (R \otimes _k K)$  and 
$\u{Q/tQ} \cong \u{N}$  in  $\sCM (R)$. 
Note that there is a natural isomorphism 
$$
\sHom _{R \otimes _k V} (M \otimes _k V, \ Q)_t \cong 
\sHom _{R \otimes _k K} (M \otimes _k K, \ Q_t). 
$$
Therefore there is a morphism $\rho :  \u{M \otimes _k V} \to \u{Q}$  in  $\sCM (R \otimes _k V)$  with  $\rho _t : \u{M \otimes _k K} \to \u{Q_t}$  is an isomorphism. 
Now take a cone of $\rho$, and we get a triangle in  $\sCM (R \otimes_kV)$; 
\begin{equation}\label{first}
\begin{CD}
\u{M\otimes _k V} @>{\rho}>>  \u{Q} @>>> \u{P'} @>>>  \u{M\otimes _kV}[1]. 
\end{CD}
\end{equation}
By the choice of  $\rho$, we have that  $\u{P'_t} \cong \u{0}$  in  $\sCM (R \otimes _k K)$, i.e. $P'_t$  is a projective $R\otimes _k K$-module. 
By virtue of Lemma \ref{Swan triangle} we have a  Cohen-Macaulay $R$-module  $X$  and a triangle in  $\sCM (R \otimes _k V)$; 
$$
\begin{CD}
\u{X \otimes _k V} @>{\mu}>> \u{X \otimes _k V} @>>> \u{P'} @>>> \u{X \otimes _k V}[1].
\end{CD}
$$ 
Utilizing the octahedron axiom, it follows from this triangle together with (\ref{first}) that there is a commutative diagram in which all rows and columns are triangles in  $\sCM (R \otimes _k V)$. 
\begin{equation}\label{tri}
\begin{CD}
@. \u{X \otimes _k V} @= \u{X\otimes _k V} \\
@. @VVV @V{\mu}VV @. \\
\u{M \otimes _k V} @>>> \u{W} @>{\nu}>> \u{X \otimes _k V} @>{\lambda}>> \u{M \otimes _k V}[1] \\
@| @VVV @VVV @| \\
\u{M \otimes _k V} @>{\rho}>> \u{Q} @>>> \u{P'} @>>> \u{M \otimes _k V}[1] \\
@. @VVV @VVV @. \\
@. \u{X \otimes _k V}[1] @= \u{X\otimes _k V}[1] \\
\end{CD}
\end{equation} 
Taking the localization by $t$, and noting that  $\u{P'_t} \cong \u{0}$, we have the following commutative diagram in which all rows and columns are triangles in  $\sCM (R \otimes _k K)$. 
$$
\begin{CD}
@. \u{X \otimes _k K} @= \u{X\otimes _k K} \\
@. @VVV @V{\mu_t}V{\cong}V @. \\
\u{M \otimes _k K} @>>> \u{W_t} @>{\nu_t}>> \u{X \otimes _k K} @>{\lambda_t}>> \u{M \otimes _k K}[1] \\
@| @VVV @VVV @| \\
\u{M \otimes _k K} @>{\rho _t}>{\cong}> \u{Q_t} @>>> \u{0} @>>> \u{M \otimes _k K}[1] \\
@. @VVV @VVV @. \\
@. \u{X \otimes _k K}[1] @= \u{X\otimes _k K}[1] \\
\end{CD}
$$
From this diagram we see that  $\nu_t$  is a splitting epimorphism, and hence $\lambda _t = 0$  in  $\sCM (R \otimes _k K)$. 
Notice that  $\lambda$  is an element of  $\Ext _{R\otimes _k V}^1 (X \otimes _k V, \ M \otimes _k V)$, and that there is a natural isomorphism 
$$
\Ext _{R\otimes _k V}^1 (X \otimes _k V, \ M \otimes _k V) \cong 
\Ext _{R}^1 (X,  M) \otimes _k V. 
$$
Thus that $\lambda _t = 0$  forces that  $t^n \lambda = 0$  in  $\Ext _{R}^1 (X, \ M) \otimes _k V$  for some  $n >0$. 
However,  since  $t$  is a non-zero divisor on  $\Ext _{R}^1 (X,  M) \otimes _k V$, this implies that  $\lambda = 0$  as an element  $\Ext _{R\otimes _k V}^1 (X \otimes _k V, \ M \otimes _k V)$. 

Now getting back to the diagram (\ref{tri}), we conclude from  $\lambda =0$  that the second row splits and that  $\u{W}$  is isomorphic to 
$\u{M \otimes _k V} \oplus \u{X \otimes _k V}$. 
Thus we have a triangle in  $\sCM (R \otimes _k V)$: 
$$
\begin{CD}
\u{X \otimes _k V} @>>> \u{M \otimes _k V} \oplus \u{X\otimes _k V}@>>> \u{Q} @>>>  \u{X \otimes _k V}[1]  \\
\end{CD}
$$
Send this triangle by the functor  $\mathcal{R} : \sCM (R \otimes _k V) \to \sCM (R)$, and we get a triangle in $\sCM (R)$  of the following form: 
$$
\begin{CD}
\u{X} @>{\binom{\u{\phi}}{\u{\psi}}}>> \u{M} \oplus \u{X}@>>> \u{N} @>>>  \u{X}[1]  \\
\end{CD}
$$
We should note that we did not use so far the assumption that $R$  is an isolated singularity.

It remains to prove that we can take a nilpotent endomorphism as $\u{\psi}$. 
For this, we apply Lemma \ref{Fitting lemma for isolated sing} to the $\u{\psi}$ above. 
As in the lemma, we have a decomposition $\u{X} = \u{X_1} \oplus \u{X_2}$  and automorphisms  $\u{\alpha}, \u{\beta}$ of  $\u{X}$ such that 
$$
\u{\alpha} \cdot \u{\psi} \cdot \u{\beta} =  \begin{pmatrix} \u{\psi_1} & 0 \\ 0 & \u{\psi_2} \\ \end{pmatrix},
$$
where  
$\u{\psi_1} : \u{X_1} \to \u{X_1}$  is an automorphism and  
$\u{\psi_2} : \u{X_2} \to \u{X_2}$  is a nilpotent endomorphism. 
Then we have an isomorphism of triangles in  $\sCM (R)$; 
$$
\begin{CD}
\u{X_1}\oplus \u{X_2} @>{\binom{\u{\phi}}{\u{\psi}}}>> \u{M} \oplus \u{X_1}\oplus \u{X_2}  @>>> \u{N} @>>>  (\u{X_1}\oplus \u{X_2})[1]  \\
@V{\u{\beta} ^{-1}}V{\cong}V    
@V{\tiny{\begin{pmatrix} 1 & 0 \\ 0 & \u{\alpha}  \\ \end{pmatrix}}}V{\cong}V    @VV{\cong}V    @VV{\cong}V     \\
\u{X_1}\oplus \u{X_2} @>>{
\tiny{
\begin{pmatrix} \u{\phi}\u{\beta}  \\ \cdots \\  \u{\psi_1} \  \  0 \\ \ \ 0 \ \ \  \u{\psi_2} \\ \end{pmatrix}
}
}> \u{M} \oplus \u{X_1}\oplus \u{X_2} @>>> \u{N} @>>>  (\u{X_1}\oplus \u{X_2})[1]  \\
\end{CD}
$$
Since  $\u{\psi _1}$  is an isomorphism, we can split $\u{X_1}$ off from the triangle in the second row above, and we get the triangle of the form; 
$$
\begin{CD}
\u{X_2} @>{\binom{\u{\phi'}}{\u{\psi_2}}}>> \u{M} \oplus \u{X_2}  @>>> \u{N} @>>>  \u{X_2}[1].
\end{CD}
$$
Since $\u{\psi _2}$  is nilpotent, this is the triangle we wanted. 
\qed

\medskip

As a direct consequence of Theorem \ref{isolated singularity}, we have the following corollary.
\begin{cor}\label{cor isolated singularity}
Let  $(R_1, \m_1, k)$ and  $(R_2, \m_2, k)$ be Gorenstein complete local $k$-algebras.
Assume that the both $R_1$ and  $R_2$  are isolated singularities, and that $k$  is an infinite field. 
Suppose there is a $k$-linear equivalence  $F : \sCM (R_1) \to \sCM (R_2)$  of triangulated categories. 
Then, for $\u{M},\  \u{N} \in \sCM (R_1)$, 
$\u{M}$  stably degenerates to $\u{N}$ if and only if  
$F(\u{M})$  stably degenerates to $F(\u{N})$. 
\end{cor}

\begin{pf}
Assume that $\u{M}$  stably degenerates to $\u{N}$ for  $\u{M}, \u{N} \in \sCM (R_1)$. 
By Theorem \ref{isolated singularity} there is a triangle 
$$
\begin{CD}
\u{Z} @>{\binom{\u{\phi}}{\u{\psi}}}>>  \u{M} \oplus \u{Z} @>>> \u{N} @>>>  \u{Z}[1], 
\end{CD}
$$
where $\u{\psi}$  is  a nilpotent element of  $\u{\rm End} _R (Z)$. 
Applying the functor $F$ to this, we have a triangle in $\sCM (R_2)$;  
$$
\begin{CD}
F(\u{Z}) @>{\binom{F(\u{\phi})}{F(\u{\psi})}}>>  F(\u{M}) \oplus F(\u{Z}) @>>> F(\u{N}) @>>>  F(\u{Z})[1]. 
\end{CD}
$$
Since  $F(\u{\psi})$ is nilpotent as well, Theorem \ref{implication} forces that $F(\u{M})$  stably degenerates to $F(\u{N})$. 
\qed
\end{pf}

\begin{rem}
Let  $(R_1, \m_1, k)$ and  $(R_2, \m_2, k)$ be Gorenstein complete local $k$-algebras as above.
Then it hardly occurs that there is a $k$-linear equivalence of categories  between  $\CM (R_1)$  and  $\CM (R_2)$. 
In fact, if it occurs, then  $R_1$  is isomorphic to  $R_2$  as a $k$-algebra. 
(See \cite[Proposition 5.1]{HY}.)
 
On the other hand, an equivalence  between $\sCM (R_1)$  and  $\sCM (R_2)$  may happen for non-isomorphic $k$-algebras. 
For example, let  $R_1 = k[[ x, y, z]]/(x^n + y^2 + z^2)$  and 
  $R_2 = k[[ x]]/(x^n)$  with characteristic of $k$ being odd  and  $n \in \N$. 
Then, by Knoerrer's periodicity (\cite[Theorem 12.10]{Y1}), we have an equivalence  
 $\sCM (k[[ x, y, z]]/(x^n + y^2 + z^2)) \cong \sCM (k[[ x]]/(x^n))$. 
Since  $k[[x]]/(x^n)$  is an artinian Gorenstein ring, the stable degeneration of modules over  $k[[x]]/(x^n)$  is equivalent to a degeneration up to free summands by Theorem \ref{implication}. 
Moreover the degeneration problem for modules over  $k[[x]]/(x^n)$  is known to be equivalent to the degeneration problem for Jordan canonical forms of square matrices of size $n$. 
Thus by virtue of Corollary \ref{cor isolated singularity}, it is easy to describe the stable degenerations of  Cohen-Macaulay modules over 
 $k[[ x, y, z]]/(x^n + y^2 + z^2)$. 
\end{rem}

%

\section{Degeneration and stable degeneration}

In this section we shall show that a stable degeneration implies a degeneration after adding some identical Cohen-Macaulay module. 
More precisely the main result of this section is the following.

\begin{thm}\label{oplusX}
Let  $R$  be a Gorenstein complete local $k$-algebra, where we assume that  $k$  is an infinite field. 
Let  $M, N \in \CM (R)$, and suppose that $\u{M}$  stably degenerates to  $\u{N}$. 
Then, there exists an  $X \in \CM (R)$  such that  $M \oplus R^m \oplus X$  degenerates to $N \oplus R^n \oplus X$ for some  $m, n \in \N$. 
\end{thm}

Before giving a proof of this theorem we need a lemma and a proposition which are basically proved by using Swan's Lemma \ref{Swan}. 
\begin{lemma}\label{Swan proective}
Let  $R$  be a noetherian local $k$-algebra, where  $k$  is an infinite field, and let  $K=k(t)$  as before. 
Given a finitely generated projective  $R \otimes _k K$-module  $P$, there are a finitely generated  $R$-module  $X$  and a non-negative integer  $n$  such that $$
P \oplus (X \otimes _k K)  \cong  (R ^n \otimes _k K) \oplus (X \otimes _k K), 
$$
as  $R\otimes _k K$-modules. 
\end{lemma}

\begin{pf}
Recall that  $R \otimes _k K = (R\otimes _k V)_t$. 
Since  $P$  is a finitely generated  $(R\otimes _k V)_t$-module, we find a finitely generated $R\otimes _k V$-submodule  $P'$  of  $P$  satisfying  $P'_t \cong  P$. 
Since  $P$ is projective, we may assume $P \subseteqq (R\otimes _k K) ^r$ for some $r \geqq 0$. 
Replacing  $P'$  with  $P' \cap (R\otimes _k V)^r$ if necessary, we may take  a submodule of a free $R\otimes _k V$-module as such a $P'$. 
Thus we can apply Lemma \ref{Swan lemma} to  $P'$ and we have an exact sequence of $R\otimes _k V$-modules;
$$
\begin{CD}  
0 @>>> X \otimes _kV @>>> (X \otimes _k V)  \oplus  (R ^n \otimes _k V) @>>> P' @>>> 0, 
\end{CD}
$$
where   $X$  is a finitely generated  $R$-module and  $n$ is a non-negative integer. 
Taking the localization by $t$  and noting that  $P = P'_t$  is projective, we have an isomorphism   
$$
P \oplus (X \otimes _k K) \cong (R^n \otimes _k K ) \oplus (X \otimes _kK).
$$ 
\qed
\end{pf}

\begin{prop}\label{Swan projective CM}
Let  $R$  be a Cohen-Macaulay local $k$-algebra, where  $k$  is an infinite field, and let  $K=k(t)$  as before. 
Given a finitely generated projective  $R \otimes _k K$-module  $P$, there are a Cohen-Macaulay $R$-module  $X$  and a non-negative integer  $n$  such that 
\begin{equation}\label{iso}
P \oplus (X \otimes _k K)  \cong  (R ^n \otimes _k K) \oplus (X \otimes _k K), 
\end{equation} 
as  $R\otimes _k K$-modules. 
\end{prop}

\begin{pf}
We have already shown that there is a finitely generated  $R$-module  $X$  which satisfies the isomorphism (\ref{iso}). 
Take an integer $\ell$  so that  $2 \ell > \dim R$ and we consider the $2 \ell$th syzygy module of  $X$, i.e. there is an exact sequence 
$$
\begin{CD}
0 @>>> \Omega ^{2\ell}X @>>>  F_{2\ell -1} @>>> \cdots @>>> F_0 @>>>  X @>>>  0, \end{CD}
$$
where each  $F_i$  is a free $R$-module.  
Since  $P$  is a projective $R \otimes _kK$-module, it is a direct summand of a free $R\otimes _kK$-module  $R^m \otimes _kK$. 
Therefore there is an idempotent  $\epsilon \in \End _{R\otimes _kK}(R^m\otimes _kK)$  such that  $P = \coker (\epsilon) = \image (1-\epsilon) = \ker (\epsilon)$. 
Thus there is an exact sequence with $2\ell +2$ terms 
$$
0 \to P \to  R^m\otimes _kK \overset{\epsilon}\to  R^m \otimes _kK \overset{1-\epsilon}\to  \cdots \overset{1-\epsilon}\to  R^m \otimes _kK \overset{\epsilon}\to R^m \otimes _kK \to P \to 0. 
$$
Tensoring the first exact sequence with $K$ over $k$ and taking the direct sum of these exact sequences, we obtain an exact sequence of the form 
$$
0 \to P \oplus (\Omega ^{2\ell}X \otimes _kK) \to G_{2\ell -1} \to \cdots \to G_0 \to P \oplus (X \otimes _kK) \to 0, 
$$
where  $G_i = R^m \otimes _kK \oplus F_i \otimes _kK$ that is a free $R \otimes _kK$-module. 
On the other hand there is an exact sequence of the form 
$$
0 \to (\Omega ^{2\ell}X \otimes _kK) \to H_{2\ell -1} \to \cdots \to H_0 \to (R^n \otimes _kK) \oplus (X \otimes _kK) \to 0, 
$$
where each $H_i$  is a free $R \otimes _k K$-module. 
Since  $(R^n \otimes _kK) \oplus (X \otimes _kK) \cong P \oplus (X \otimes _kK)$, it follows from Schanuel's lemma that 
there is an isomorphism of the form 
$$
P \oplus (\Omega ^{2\ell}X \otimes _kK) \oplus (R^a \otimes _kK) 
\cong (\Omega ^{2\ell}X \otimes _kK) \oplus (R^b \otimes _kK), 
$$
for some $a, b \in \N$. 
Note that in such a situation we must have  $a < b$ whenever  $P \not= 0$. 
Thus setting  $X' = \Omega ^{2\ell }X \oplus R^a$, we have 
$$
P \oplus (X' \otimes _kK ) \cong (R^{b-a} \otimes _kK) \oplus (X' \otimes _kK)
$$
and clearly  $X'$  is a Cohen-Macaulay module over  $R$. 
\qed
\end{pf}

\vspace{12pt}

Now we can prove Theorem \ref{oplusX}. 

\vspace{12pt}
[Proof of Theorem \ref{oplusX}]

We assume that  $\u{M}$ stably degenerates to $\u{N}$  for  $M, N \in \CM (R)$. 
By definition there is a  Cohen-Macaulay $R\otimes _k V$-module  $Q$  such that  $\u{Q_t} \cong \u{M \otimes _k K}$  in $\sCM (R\otimes _k K)$  and  $\u{Q/tQ} \cong \u{N}$ in  $\sCM (R)$. 
Thus  $Q_t \oplus P_1 \cong (M \otimes _k K) \oplus P_2$  in $\CM (R \otimes _kK)$  for projective  $R \otimes _k K$-modules $P_1$  and  $P_2$, and 
$Q/tQ \oplus R^a \cong N \oplus R^b$  in  $\CM (R)$  for some $a, b \in \N$. 
It then follows from Proposition \ref{Swan projective CM} there are  Cohen-Macaulay $R$-modules  $X_1, X_2$  and integers  $n_1, n_2$  satisfying 
$$
P_i  \oplus (X_i \otimes _k K)  \cong  (R ^{n_i} \otimes _k K) \oplus (X_i \otimes _k K), 
$$
for $i = 1,2$. 
Now setting  $X = X_1 \oplus X_2$ we have an isomorphism in $\CM ( R \otimes _k K)$; 
$$
Q_t \oplus (R ^{n_1}\otimes _k K)  \oplus (X \otimes _kK) \cong 
(M \otimes _k K)  \oplus (R ^{n_2}\otimes _k K)  \oplus (X \otimes _kK). 
$$ 
We denote by  $\tilde{Q}$ the  Cohen-Macaulay module  
$Q \oplus (R^{a+n_1}\otimes _k V)  \oplus (X \otimes _k V)$ over  $R \otimes _k V$. 
Then it follows 
$$
\tilde{Q}_t \cong (M \otimes _k K)  \oplus (R ^{a+n_2}\otimes _k K)  \oplus (X \otimes _kK), 
$$
in  $\CM (R \otimes _kK)$  and $\tilde{Q}/t \tilde{Q} \cong N \oplus R ^{b+n_1}\oplus X$  in  $\CM (R)$.  
Therefore $M \oplus R ^{a+n_2} \oplus X$  degenerates to $N \oplus R ^{b+n_1}\oplus X$.
\qed  

\medskip
The converse of Theorem \ref{oplusX} does not hold even in the case when  $R$  is artinian. 
The following example is taken from the paper \cite{R} of Riedtmann. 

\begin{ex}\label{Riedtmann}
Let  $R = k[[x, y]]/(x^2, y^2)$. 
Note that  $R$  is an artinian Gorenstein local ring. 
Now consider the modules 
$M_{\lambda} = R/(x - \lambda y)R$  for all $\lambda \in k$. 
We denote by $k$ the unique simple module  $R/(x,y)R$ over  $R$. 
In this case, it is known by \cite[Example 3.1]{R} that 
$R \oplus k^2$  degenerates to  $M_{\lambda} \oplus M_{\mu} \oplus k^2$  for any choice of  $\lambda , \mu \in k$. 
 
We claim that  $\u{R}$  never stably degenerates to  $\u{M_{\lambda} \oplus M_{\mu}}$  if  $\lambda + \mu \not= 0$. 

In fact, if there is such a stable degeneration, then it follows from Theorem \ref{implication} that  $R^m$  degenerates to  $M_{\lambda} \oplus M_{\mu} \oplus R^n$  for some  $m, n \in \N$. 
Since $[R^m] = [M_{\lambda} \oplus M_{\mu} \oplus R^n]$ in the Grothendieck group, we have  $m > n \geqq 0$. 
Now we apply Theorem \ref{decreaseFitting} to obtain an inclusion of Fitting ideals; 
$\F_n^R(M_{\lambda} \oplus M_{\mu} \oplus R^n) \subseteq \F_n^R (R^m)$. 
We note that  $\F_n^R(R^m) = 0$ since  $n <m$, and an easy computation shows that 
$$
\F_n^R(M_{\lambda} \oplus M_{\mu} \oplus R^n) = 
\F_0^R(M_{\lambda})\F_0^R(M_{\mu})\F_n^R(R^n) = 
(x-\lambda y)(x-\mu y)R = (\lambda + \mu)xyR. 
$$
Hence we must have $\lambda + \mu =0$. 
\end{ex}

\section{Stable degeneration order}

As an application of Theorem \ref{oplusX} we can define the stable degeneration order for Cohen-Macaulay modules.

\begin{defn}\label{def st deg order}
Let  $(R, \m, k)$  be a Gorenstein complete local $k$-algebra as before, and let $\u{M}, \ \u{N} \in \sCM (R)$. 
If there is a sequence of objects  $\u{L_0}, \u{L_1}, \u{L_2}, \ldots, \u{L_n}$  in  $\sCM (R)$  such that 
$\u{L_0}=\u{M}$, $\u{L_n}=\u{N}$  and 
$\u{L_i}$  stably degenerates to  $\u{L_{i+1}}$  for  $i=0,1,\cdots ,n-1$, 
then 
we write  $\u{M} \leq_{st} \u{N}$. 
\end{defn}

Theorem \ref{oplusX} shows that the relation $\leq_{st}$  gives a partial order on the set of isomorphism classes of objects in  $\sCM (R)$. 
In fact we can prove the antisymmetric law  for  $\leq_{st}$.

\begin{thm}\label{st deg order}
Let  $(R, \m, k)$  be a Gorenstein complete local algebra over an infinite field $k$, and let $\u{M}, \ \u{N} \in \sCM (R)$. 
If $\u{M} \leq_{st} \u{N}$  and  $\u{N} \leq_{st} \u{M}$, then  $\u{M} \cong \u{N}$. 
\end{thm}

\begin{pf}
Since  $\u{M} \leq _{st} \u{N}$, there are $\u{L_0}, \u{L_1}, \u{L_2}, \ldots, \u{L_n}$  in  $\sCM (R)$  such that $\u{L_0}=\u{M}$, $\u{L_n}=\u{N}$  and 
$\u{L_i}$  stably degenerates to  $\u{L_{i+1}}$  for  $i=0,1,\cdots ,n-1$. 
It follows from Theorem \ref{oplusX} that 
$L_i \oplus R^{a_i} \oplus X_i$  degenerates to  
$L_{i+1} \oplus R^{b_{i}} \oplus X_i$, where $X_i \in \CM(R)$  and  $a_i, \ b_i \in \N$  for each   $i=0,1,\cdots ,n-1$. 
Then, set  
$$
L'_i  = L_i \oplus R^{b_0+\cdots + b_{i-1} + a_i + \cdots +a_{n-1}} \oplus X_0 \oplus \cdots \oplus X_{n-1}, 
$$
and we can see that  
$L'_i$ degenerates to  $L'_{i+1}$  for   $i=0,1,\cdots ,n-1$. 
Therefore, we have  $L'_0 \leq _{deg} L'_n$  under the ordinary degeneration order, thus setting  $a = a_0+\cdots +a_{n-1}$, $b = b_0+\cdots +b_{n-1}$  and  $X = X_1 \oplus \cdots \oplus X_{n-1}$,  we have 
$$
M \oplus R^{a} \oplus X \leq_{deg} N \oplus R^{b} \oplus X. 
$$
(See \cite[Definition 4.11]{Y2} for the detail of degeneration order $\leq_{deg}$.) 
Similarly, using the assumption that  $\u{N} \leq _{st} \u{M}$, we get 
$$
N \oplus R^{c} \oplus Y \leq_{deg} 
M \oplus R^{d} \oplus Y, 
$$
for some  $c, d \in \N$  and  $Y \in \CM (R)$. 
Thus we have the inequality 
$$
M \oplus  R^{a+c} \oplus X \oplus Y
 \leq_{deg} 
N \oplus R^{b+c}  \oplus X \oplus Y
 \leq_{deg} 
M \oplus  R^{b+d} \oplus X \oplus Y. 
$$
Since there is a degeneration, we have the equality 
$[M \oplus  R^{a+c} \oplus X \oplus Y] = [M \oplus  R^{b+d} \oplus X \oplus Y
]$  in the Grothendieck group  $K_0(\mod (R))$. 
See Remark \ref{if degerates}.
Thus it follows that $[R^{a+c}] = [R^{b+d}]$, hence $a+c = b+d$. 
Since $\leq_{deg}$  satisfies the antisymmetric law (\cite[Theorem 2.2, Proposition 4.4]{Y2}), 
we have the isomorphism 
$$
M \oplus  R^{a+c} \oplus X \oplus Y  \cong N \oplus R^{b+c}  \oplus X \oplus Y. $$
Note that  $\CM (R)$  is a Krull-Schmidt category, and thus this isomorphism forces  
$M \oplus  R^{a} \cong N \oplus R^{b}$. 
Therefore  $\u{M} \cong \u{N}$  in $\sCM (R)$. 
\qed
\end{pf}


If $R$ is an isolated singularity, then we can prove that the stable degeneration for Cohen-Macaulay modules is transitive. 
More precisely we can prove the following proposition.

\begin{prop}
Let  $(R, \m, k)$  be a Gorenstein complete local algebra over an infinite field $k$, and assume that  $R$  has only an isolated singularity. 

For $\u{L}, \u{M}, \u{N} \in \sCM (R)$, 
if $\u{L}$ stably degenerates to $\u{M}$ and if $\u{M}$ stably degenerates to $\u{N}$, then  $\u{L}$ stably degenerates to $\u{N}$. 
In particular, the partial order $\u{M} \leq_{st} \u{N}$  is equivalent to saying that  $\u{M}$ stably degenerates to $\u{N}$. 
\end{prop}

\begin{pf}
The proof proceeds as in the same way as the proof of \cite[Theorem 2.1]{Z2}. 
In fact, if there are triangles 
$$
\begin{CD}
\u{Z_1} @>{\binom{\u{\phi_1}}{\u{\psi_1}}}>>  \u{L} \oplus \u{Z_1} @>>> \u{M} @>>>  \u{Z_1}[1], \\
\u{Z_2} @>{\binom{\u{\phi_2}}{\u{\psi_2}}}>>  \u{M} \oplus \u{Z_2} @>>> \u{N} @>>>  \u{Z_2}[1], 
\end{CD} 
$$
with  $\u{\psi_1}$ and $\u{\psi_2}$ being nilpotent, then we can construct a new triangle of the form;
$$
\begin{CD}
\u{Z_3} @>{\binom{\u{\phi_3}}{\u{\psi_3}}}>>  \u{L} \oplus \u{Z_3} @>>> \u{N} @>>>  \u{Z_3}[1],
\end{CD} 
$$
where  $\u{\psi_3}$  may not be nilpotent. 
(One can prove this completely in a similar way to the proof of \cite[Theorem 2.1]{Z2} but by replacing short exact sequences there by triangles.) 
As in the same manner of the last half of the proof of Theorem \ref{isolated singularity}, we can replace $\u{\psi_3}$ by a nilpotent endomorphism by utilizing Lemma \ref{Fitting lemma for isolated sing}. 
\qed
\end{pf}

\bibliographystyle{amsplain}

\ifx\undefined\bysame
\newcommand{\bysame}{\leavevmode\hbox to3em{\hrulefill}\,}
\fi


\end{document}